\documentclass[12pt]{article}
\usepackage{amsmath}
\usepackage{amsfonts}
\usepackage{amssymb}
\usepackage{theorem}

\title{Continuous structures of quantum circuits}
\author{A. Ivanov
\thanks{The author is supported by Polish National Science Centre grant DEC2011/01/B/ST1/01406} } 
\date{ } 

\newtheorem{thm}{Theorem}[section] 
\newtheorem{lem}[thm]{Lemma}
\newtheorem{definicja}[thm]{Definition}

\newtheorem{remark}[thm]{Remark}

\begin{document}
\maketitle
\topskip 20pt

\begin{quote}
{\bf Abstract.} 
We consider continuous structures which are obtained 
from finite dimensional Hilbert spaces over $\mathbb{C}$ 
by adding some unitary operators. 
Quantum automata and circuits are naturally interpretable 
in such structures. 
We consider appropriate algorithmic problems concerning 
continuous theories of natural classes of these structures. 

\bigskip


\end{quote}

\bigskip


\section{Introduction} 

Continuous logic has become the basic model theoretic 
tool for Hilbert spaces and $\mathbb{C}^*$-algebras: 
see \cite{BYBHU}, \cite{BYU} and \cite{FHS}. 
This suggests that quantum circuits, quantum automata 
and quantum computations in general 
can be defined in appropriate continuous structures 
and studied by means of continuous logic.  
The paper is an attempt of this approach. 
\parskip0pt 

It is worth noting that typical continuous structures 
appearing in quantum informatics are finite dimensional, 
i.e. compact and not very interesting from the point of view 
of model theory.  
On the other hand we demonstrate in our paper 
that study of continuous theories of {\em classes} of 
these structures may be promising.    
\parskip0pt 

We remind the reader that states of quantum systems 
are represented by normed vectors of tensor products 
$$ 
(...(\mathcal{B}_1 \bigotimes \mathcal{B}_2 )\bigotimes .... )\bigotimes \mathcal{B}_k ,
$$ 
$$
\mbox{  where } \mathcal{B}_i \cong \mathbb{C}\bigoplus \mathbb{C} 
\mbox{ under isomorphisms of Hilbert spaces, } i\le k. 
$$ 
In Dirac's notation elements of $\mathcal{B}_i$ 
are denoted by $|h\rangle$ and 
tensors 
$$ 
(...(|h_1 \rangle \otimes |h_2 \rangle ) ...) \otimes |h_k \rangle 
\mbox{ are denoted by } 
|h_1 h_2 ...h_k \rangle . 
$$      
Any normed $h\in \mathcal{B}_i$ is called a {\bf qubit};  
it is a linear combination of  
$|0\rangle = (1,0)$ and $|1\rangle =(0,1)$.  

The probability amplitude $a (\phi \rightarrow \psi)$ 
is defined as the inner product 
$\langle \psi |\phi\rangle$   
and the probability $p(\phi \rightarrow \psi)$ 
is $|a(\phi \rightarrow \psi)|^2$. 
Dynamical evolutions of the quantum system are 
represented by unitary operators on  
$\mathcal{B}^{\otimes k}$. 

We call the structure $\mathcal{B}^{\otimes k}$ 
enriched by unitary operators $U_1 ,...,U_t$ 
a {\bf dynamical $n$-qubit space}.  
It in particular defines a family of quantum automata  
over the language 
$\{ 1,...,t\}^{*}$, where each automaton 
is determined by the $2^n$-dimensional  
diagonal matrix $P$ of the projection to final states. 
Fixing  $\lambda\in \mathbb{Q}$ it is said that 
a word $w=i_1 ...i_k$ is {\em accepted} if 
$$
\parallel P U_{i_k} ... U_{i_1} |0^{\otimes n}\rangle \parallel^2 >\lambda . 
$$   
These issues are described in \cite{gruska}, \cite{MC} and \cite{DJK}. 

We will consider dynamical $n$-qubit spaces in 
continuous logic. 
All necessary information on  
continuous logic  will be described in 
the next section. 
The main results of the paper (see Section 4) 
concern decidability of continuous theories of 
classes of dynamical qubit spaces. 
They are  in particular motivated by \cite{DJK} 
(see Section 2 below). 
Section 3 contains some general 
observations concerning decidability.  
We think that this section is interesting by itself. 
It is naturally connected with the material 
of \cite{BYP}, \cite{DGP} and \cite{GH}.

It is worth noting that  continuous logic can be considered as 
a theory in some extension ( {\bf RPL}$\forall$ ) 
of {\L}ukasiewicz logic (see \cite{DGP}). 
The latter is traditionally linked with quantum mechanics, 
\cite{CCGP}, \cite{pykacz}. 
Thus the idea that continuous logic should enter 
into the field is quite natural. 
On the other hand the author thinks that the context 
of the paper is original.

\section{Continuous structures.} 

We fix a countable continuous signature 
$$
L=\{ d,R_1 ,...,R_k ,..., F_1 ,..., F_l ,...\}. 
$$  
Let us recall that a {\em metric $L$-structure} 
is a complete metric space $(M,d)$ with $d$ bounded by 1, 
along with a family of uniformly continuous operations on $M$ 
and a family of predicates $R_i$, i.e. uniformly continuous maps 
from appropriate $M^{k_i}$ to $[0,1]$.   
It is usually assumed that to a predicate symbol $R_i$ 
a continuity modulus $\gamma_i$ is assigned so that when 
$d(x_j ,x'_j ) <\gamma_i (\varepsilon )$ with $1\le j\le k_i$ 
the corresponding predicate of $M$ satisfies 
$$ 
|R_i (x_1 ,...,x_j ,...,x_{k_i}) - R_i (x_1 ,...,x'_j ,...,x_{k_i})| < \varepsilon . 
$$ 
It happens very often that $\gamma_i$ coincides with $id$. 
In this case we do not mention the appropriate modulus. 
We also fix continuity moduli for functional symbols. 
Note that each countable structure can be considered 
as a complete metric structure with the discrete $\{ 0,1\}$-metric. 

By completeness continuous substructures of a continuous structure are always closed subsets. 

Atomic formulas are the expressions of the form $R_i (t_1 ,...,t_r )$, 
$d(t_1 ,t_2 )$, where $t_i$ are terms (built from functional $L$-symbols). 
In metric structures they can take any value from $[0,1]$.   
{\em Statements} concerning metric structures are usually 
formulated in the form 
$$
\phi = 0 
$$ 
(called an $L$-{\em condition}), where $\phi$ is a {\em formula}, 
i.e. an expression built from 
0,1 and atomic formulas by applications of the following functions: 
$$ 
x/2  \mbox{ , } x\dot- y= max (x-y,0) \mbox{ , } min(x ,y )  \mbox{ , } max(x ,y )
\mbox{ , } |x-y| \mbox{ , } 
$$ 
$$ 
\neg (x) =1-x \mbox{ , } x\dot+ y= min(x+y, 1) \mbox{ , } x \cdot y \mbox{ , } sup_x \mbox{ and } inf_x . 
$$ 
A {\em theory} is a set of $L$-conditions without free variables 
(here $sup_x$ and $inf_x$ play the role of quantifiers). 
   
It is worth noting that any formula is a $\gamma$-uniformly continuous 
function from the appropriate power of $M$ to $[0,1]$, 
where $\gamma$ is the minimum of continuity moduli of $L$-symbols 
appearing in the formula. 

The condition that the metric is bounded by $1$ is not necessary. 
It is often assumed that $d$ is bounded by some rational number $d_0$. 
In this case the (dotted) functions above are appropriately modified.  
Sometimes predicates of continuous structures map $M^n$ to some 
$[q_1 ,q_2 ]$ where $q_1 ,q_2 \in \mathbb{Q}$.

\begin{remark} \label{FaHaSh} 
{\em 
Following Section 4.2 of \cite{FHS} we define a topology 
on $L$-formulas relative to a given continuous theory $T$. 
For $n$-ary formulas $\phi$ and $\psi$ of the same sort set 
$$ 
{\bf d}^T_{\bar{x}} (\phi ,\psi) = sup \{ |\phi (\bar{a}) -\psi (\bar{a} )|: \bar{a} \in M, M\models T\} . 
$$ 
The function ${\bf d}^T_{\bar{x}}$ is a pseudometric. 
The language $L$ is called {\em separable} if 
for every $L$-theory $T$ and any tuple $\bar{x}$ 
the density character of ${\bf d}^T_{\bar{x}}$ is countable. 
By Proposition 4.5 of \cite{FHS} in this case for every 
$L$-model $M$ the set of all interpretations of $L$-formulas 
in $M$ is separable in the uniform topology. 
}
\end{remark}

The paper \cite{BYP} gives  
fourteen axioms of continuous 
first order logic, denoted by  
(A1) - (A14), and the correspobding version of 
{\em modus ponens}: 
$$ 
\frac{\phi \mbox{ , } \psi \dot{-} \phi}{\psi} . 
\mbox{ where } \phi , \psi \mbox{ are continuous formulas.} 
$$ 
Corollary 9.6 of \cite{BYP} states 
\begin{quote} 
Let $\Gamma$ be a set of continuous formulas 
of a continuous signature $L$ with a metric. 
Let $\phi$ be a continuous $L$-formula. 
Then the following conditions are equivalent: \\ 
(i) for any continuous structure $M$ and any $M$-assignment of variables, 
if $M$ satisfies all statements $\psi = 0$, 
$\psi \in \Gamma$, then $M$ satisfies $\phi =0$; \\ 
(ii) $\Gamma \vdash \phi \dot{-} 2^{-n}$ for all $n\in \omega$. 
\end{quote} 
It is called {\bf approximated strong 
completeness for continuous first-order logic}. 
The following statement is Corollary 9.8 from \cite{BYP}. 
\begin{quote} 
Under circumstances above 
the following values are the same: \\  
(i) $sup \{ \phi^{M} :  \mbox{ for all } M \models \Gamma =0\}$; \\ 
(ii) $inf \{ p\in \mathbb{Q}: \Gamma \vdash   \phi \dot{-} p\}$. 
\end{quote} 
We denote this value by $\phi^{\circ}$.  

If the language $L$ is computable, the set of 
all continuous $L$-formulas and the set of all $L$-conditions of the form 
$$ 
\phi \le \frac{m}{n} \mbox{ , where } \frac{m}{n}\in \mathbb{Q}_+ , 
$$ 
are computable. 
Moreover if $\Gamma$ is a computably emunerable 
set of formulas, then the relation $\Gamma \vdash \phi$ 
is computably enumerable. 

Corollary 9.11 of \cite{BYP} states that when $\Gamma$ 
is computably enumerable 
and $\Gamma=0$ axiomatises a complete theory, 
then the value of $\phi$ with respect to $\Gamma$ is a recursive real
which is  uniformly computable from $\phi$.  
This exactly means that the corresponding 
complete theory is {\em decidable} (see Section 3). 
Note that in this case the value of $\phi$ coincides with 
$\phi^{\circ}$.

\section{Dynamical $n$-qubit spaces} 

We treat a Hilbert space over $\mathbb{R}$ 
exactly as in Section 15 of \cite{BYBHU}. 
We identify it with a many-sorted metric structure 
$$
(\{ B_n\}_{n\in \omega} ,0,\{ I_{mn} \}_{m<n} ,
\{ \lambda_r \}_{r\in\mathbb{R}}, +,-,\langle \rangle ),
$$
where $B_n$ is the ball of elements of norm $\le n$, 
$I_{mn}: B_m\rightarrow B_n$ is the inclusion map, 
$\lambda_{r}: B_m\rightarrow B_{km}$ is scalar 
multiplication by $r$, with $k$ the unique integer 
satisfying $k\ge 1$ and $k-1 \le |r|<k$; 
futhermore, $+,- : B_n \times B_n \rightarrow B_{2n}$ 
are vector addition and subtraction and 
$\langle \rangle : B_n \rightarrow [-n^2 ,n^2 ]$ 
is the predicate of the inner product. 
The metric on each sort is given by 
$d(x,y) =\sqrt{ \langle x-y, x-y \rangle }$.   
For every operation the continuity modulus is standard.  
For example in the case of $\lambda_r$ this is $\frac{z}{|r|}$. 

Stating existence of infinite approximations 
of orthonormal bases by axioms of the form 
$$ 
inf_{x_1 ,...,x_n} max_{1\le i<j\le n} (|\langle x_i ,x_j\rangle -\delta_{i,j} |) =0 \mbox{ , } n\in \omega , 
$$
$$ 
\delta_{i,j} \in \{ 0,1\} \mbox{ with } \delta_{i,j} =1 \leftrightarrow i=j , 
$$   
we axiomatise  infinite dimensional Hilbert spaces. 
By \cite{BYBHU} they form the class of models of a complete 
theory which is $\kappa$-categorical for all infinite $\kappa$, 
and admits elimination of quantifiers. 

When we assume that the space is finite dimensional 
all sorts $B_n$ become compact. 
This case corresponds to the case 
of finite structures in ordinary 
model theory. 
The dimension can be described by 
maximal $n$ so that the following sentence holds.  
$$ 
inf_{y_1 ,...,y_n } sup_x (| ( |\langle x,x\rangle |^2 - |\langle x,y_1 \rangle |^2 - ... 
...-|\langle x, y_n \rangle |^2 ) |)=0 .
$$ 
The corresponding continuous theory 
admits elimination of quantifiers. 
This follows by the argument of  
Lemma 15.1 from \cite{BYBHU}.

This approach can be naturally extended to complex Hilbert spaces, 
$$
(\{ B_n\}_{n\in \omega} ,0,\{ I_{mn} \}_{m<n} ,
\{ \lambda_c \}_{c\in\mathbb{C}}, +,-,\langle \rangle_{Re} , \langle \rangle_{Im} ). 
$$
We only extend the family 
$\lambda_{r}: B_m\rightarrow B_{km}$, $r\in \mathbb{R}$, 
to a family $\lambda_{c}: B_m\rightarrow B_{km}$, $c\in \mathbb{C}$, 
of scalar products by $c\in\mathbb{C}$, with $k$ 
the unique integer satisfying $k\ge 1$ and $k-1 \le |c|<k$. 

We also introduce $Re$- and $Im$-parts of the inner product. 

If we remove from the signature of complex Hilbert spaces 
all scalar products by $c\in \mathbb{C}\setminus \mathbb{Q}[i]$, 
we obtain a countable subsignature 
$$
(\{ B_n\}_{n\in \omega} ,0,\{ I_{mn} \}_{m<n} ,
\{ \lambda_c \}_{c\in\mathbb{Q}[i]}, +,-,\langle \rangle_{Re} , \langle \rangle_{Im} ),
$$
which is {\em dense} in the original one: \\ 
if we present $c\in \mathbb{C}$ by a sequence $\{ q_i \}$ 
from $\mathbb{Q}[i]$ converging to $c$, 
then the choice of the continuity moduli of 
the restricted signature still guarantees that 
in any sort $B_n$ the functions $\lambda_{q_i}$ 
form a sequence which converges to  $\lambda_c$ 
with respect to the metric 
$$ 
sup_{x\in B_n} \{ |f^M (x) - g^M (x)| : M \mbox{ is an $L$-structure } \}.  
$$  
This obviously implies that the original language 
of Hilbert spaces is {\em separable}. 
In particular we may apply Remark \ref{FaHaSh}. 

\bigskip 

We extend structures of complex Hilbert spaces by 
additional discrete sort $Q$ with $\{ 0,1\}$-metric  
and a map $qu: Q\rightarrow B_1$ so that 
the set $qu(Q)$ is an orthonormal basis of $\mathbb{H}$. 

When $Q$ consists of $2^n$ elements we may denote them 
by $|i_0 ...i_{n-1} \rangle$ witn $i_j \in \{ 0,1\}$. 
In fact this is the $n$-{\bf qubit space}. 
It is distinguished by by the axiom 
$$ 
sup_x (| ( |\langle x,x\rangle |^2 - |\langle x,qu(q_0 )\rangle |^2 - ... 
...-|\langle x,qu(q_i )\rangle |^2 - ... - |\langle x,qu(q_{2^n-1} )\rangle |^2) |)=0 
$$ 

\bigskip 

To study dynamical evolutions of quantum cirquits 
we introduce the following expansion of 
$n$-qubit spaces.
Let us fix a natural number $t$ and consider 
the class of {\bf dynamical $n$-qubit spaces} 
(i.e. $2^n$-dimensional) in the extended signature 
$$
(Q, qu , \{ B_n\}_{n\in \omega} ,0,\{ I_{mn} \}_{m<n} ,
\{ \lambda_c \}_{c\in\mathbb{Q}[i]}, +,-,\langle \rangle_{Re} , \langle \rangle_{Im}, U_1 ,...,U_t  ),
$$
where $U_j$, $1\le j \le t$,  are  symbols of unitary 
operators of $\mathbb{H}$.  
We may assume that all $U_j$ are defined only on $B_1$. 

\begin{lem} \label{comp}
The complete continuous theory of each structure 
of this form is axiomatised by the standard 
axioms of $n$-qubit spaces, the axioms stating that each 
$U_j$ is a unitary operator and the axioms 
describing the matrix of $U_j$ in the 
basis $qu (Q)$ (appropriately enumerated by $1,...,N=2^n$): 
$$ 
inf_{q_1 ...q_N} sup_{j,l} (\parallel U_j (qu (q_l ))-\sum \lambda_{c_k}(qu (q_k )) \parallel ) \le \varepsilon_l 
\mbox{ , where } \varepsilon_l \in \mathbb{Q}, c_k \in \mathbb{Q}[i] . 
$$ 
\end{lem} 

Indeed in any model with these axioms the values of 
$U_j (qu (q_l ))$ have the same coordinates in 
appropriately enumerated $qu(Q)$. 
Thus the lemma is obvious. 

Each dynamical $n$-qubit space defines 
a family of quantum automata  over the language 
$\{ 1,...,t\}^{*}$, where each automaton 
is determined by the $2^n$-dimensional  
diagonal matrix $P$ of the projection to final states
of $qu(Q)$. 

Fixing  $\lambda\in \mathbb{Q}$ we say that 
a word $w=i_1 ...i_k$ is {\em accepted} 
by the corresponding $P$-automaton if 
$$
ACC_w = \parallel P U_{i_k} ... U_{i_1} |0^{\otimes n}\rangle \parallel^2 >\lambda . 
$$   
The corresponding algorithmic problems 
were in particular studied in the paper of 
H.Derksen, E.Jeandel, P.Koiran  \cite{DJK}.  
They have proved that the following problems are decidable for 
$U_1 ,...,U_t$  over finite extensions of $\mathbb{Q}[i]$: \\ 
(i) Is there $w$ such that $ACC_w >\lambda$? \\ 
(ii) Is a threshold $\lambda$ isolated, i.e. is there $\varepsilon$ 
that for all $w$,  $|ACC_w - \lambda |\ge \varepsilon$ ? \\ 
(iii) Is there a threshold $\lambda$ which is isolated? 

The observation that  given $P$ each statement 
$ACC_w \le \lambda$  or $|ACC_w - \lambda |\ge \varepsilon$
can be rewritten as a continuous statement 
of the theory of dynamical $n$-qubit spaces 
partially motivated our research in this paper. 
\begin{quote}  
{\em Describe classes of qubit spaces 
(possibly with additional operators) 
having decidable continuous theory. }
\end{quote} 
Below we study continuous statements $\theta$ so that 
the continuous theory of (dynamical) qubit spaces 
satisfying $\theta$ is decidable. 

\section{Decidability/undecidability of continuous theories} 

In this section we assume that the signature $L$ 
is computable and values of formulas are in $[0,1]$. 
The interval $[0,1]$ can be obviously replaced 
by any compact interval. 
We start with the following definition from 
\cite{BYP}. 

\begin{definicja} 
A continuous theory $T$ is called {\bf decidable} 
if for every sentence $\phi$ the {\em degree of truth}  
$$ 
\phi^{\circ} = sup \{ \phi^{M} :   M \models T\}
$$ 
is a computable real which is uniformly computable from $\phi$. 
\end{definicja} 
This exactly means that there is an algorithm which 
for every $\phi$ and a rational number $\delta$ 
finds a rational $r$ such that 
$|r-\phi^{\circ}|\le \delta$. 

Note that decidability of $T$ does not imply that 
the set of all continuous $\phi$ with $\phi^M = 0$ 
for all $M\models T$, is computable. 
On the other hand it is easy to see  that deciadability 
of $T$ follows from this condition.  

The following theorem is a counterpart of Ershov's 
decidability criterion (Theorem 6.1.1 of \cite{ershov}). 
Here we call a sequence of complete continuous theories $\{T_i, i\in\omega \}$  
{\bf effective} if the relation   
$$ 
\{ (\theta ,j): \theta \mbox{ is a statement so that } T_j \vdash \theta \}
$$ 
is computably enumerable.    

\begin{thm} \label{Ersov}
A continuous theory $T$ is decidable if and only if 
$T$ can be defined by a computably enumerable system of 
axioms and $T$ can be presented $T=\bigcap_{i\in\omega} T_i$  
where $\{T_i, i\in\omega \}$ is an effective sequence of 
complete continuous theories.    
\end{thm} 
 
{\em Proof.}  Sufficiency.  
Let $\phi$ be a continuous sentence. 
For every natural $n$ we can apply an effective procedure
which looks for conditions of the form $\phi \le \frac{k}{n}$ 
derived from the axioms of $T$ and conditions of 
the form $\frac{l}{n} \le \phi$ which appear in 
some $T_j \vdash \frac{l}{n} \le \phi$. 
This always gives a number $k<n-1$ such that 
$\frac{k}{n} \le \phi \le \frac{k+2}{n}$.  

Necessity.
For every sentence $\phi$ we fix a computably enumerable 
sequence of segments $[l_{n,\phi}, r_{n,\phi}]$ converging 
to $\phi^o$ so that $\phi^{\circ} \in [l_{n,\phi}, r_{n,\phi}]$. 
Then all statements $\phi\le r_{n,\phi}$ form a computably 
enumerable sequence of axioms of $T$.

Now for every sentence $\phi$ we effectively 
build a complete theory $T_{n,\phi}\supset T$ with 
$T_{n,\phi} \vdash l_{n,\phi} \dot{-} \phi \le 2^{-n}$. 
In fact such a construction produces 
an effective family $T_i$, $i\in \omega$, 
from the formulation. 
Indeed, then for every natural $n$ 
we can find a sufficiently large $m$ so that 
$T_{m,\phi} \vdash \phi^{\circ} \dot{-} \phi \le 2^{-n}$ 
(here $\phi^{\circ}$ is defined by $T$). 
This obviously implies that 
$T$ coincides with the intersection 
of all $T_{m,\phi}$. 
Effectiveness will be verified below. 

At Step 0 we define $T_{n,\phi ,0}$ 
to be the extension of $T$ by the axiom 
$l_{n,\phi} \dot{-} \phi \le 0$.  
At every step $m+1$ we build 
a finite extension $T_{n,\phi ,m+1}$ 
of $T$ so that each inequality 
$\psi \le 0$ from $T_{n,\phi ,m} \setminus T$ 
is transformed into an inequality 
$\psi \le \varepsilon$, where 
$\varepsilon \le 2^{-(2n+m+1)}$. 
The 'limit theory' 
$T_{n, \phi} = lim_{m\rightarrow \infty} T_{n,\phi, m}$ 
is defined by the limits of 
these values $\varepsilon$  
for all formulas $\psi$.  
Note that it can happen that $\varepsilon \le 0$, 
i.e. the transformed inequality 
is of the form $\psi+\delta \le 0$, 
with $\delta>0$. 
On the other hand we will see that 
for every $\psi$ the axioms of  
$lim_{m\rightarrow \infty} T_{n,\phi, m}$ 
give an effective sequence of rational numbers 
which converges to the value of $\psi$ 
under this theory. 

Let us enumerate all triples $(n,\phi ,\psi )$ 
by natural numbers $>0$ so that each triple 
has infinitely many numbers. 
Assume that the number $m+1$ codes a triple $(n, \phi , \psi )$. 
For all $n'\not= n$ we put $T_{n', \phi' ,m+1} = T_{n', \phi' ,m}$.  
Assume that at Step $m$ the theory 
$T_{n,\phi ,m}\setminus T$ already contains 
inequalities $\frac{k_l}{l} \le \psi_l \le \frac{k'_l }{l}$ 
for some natural $l$ and $k_l , k'_l \le l$.  
In particular we admit that the $0$-th inequality 
$l_{n,\phi} \dot{-} \phi \le 0$ has been 
already transformed into an inequality 
$l_{n,\phi} \dot{-} \phi \le \varepsilon$  
for some $\varepsilon \le \sum_{i\le m} 2^{-(2n+i+1)}$. 
Let $\theta$  be 
$$
\psi \dot{-} 2^{2n+m+2} max_l ( max(\psi_l \dot{-} \frac{k'_l }{l}, 
\frac{k_l}{l} \dot{-} \psi_l ,l_{n,\phi} \dot{-} (\phi + \sum_{i\le m} 2^{-(2n+i+1)} ))). 
$$ 
Since $T$ is decidable we compute 
$k_{m+1} <m$ so that 
$\frac{k_{m+1}}{m+1} \le \theta^{\circ} \le \frac{k_{m+1}+2}{m+1}$.  
Then the value of $\psi$ under $T_{n,\phi ,m}$ 
is not greater than $\frac{k_{m+1}+2}{m+1}$. 
This means that extending $T_{n,\phi ,m}$ by 
$\psi \le \frac{k_{m+1}+2}{m+1}$ 
we preserve consistency of the theory.
If $k_{m+1} =0$ this finishes our construction at this step. 

If $k_{m+1}>0$ we need an additional 
correction. 
Let $\theta'$  be 
$$
\psi \dot{-} 2^{2n+m+2} max_l ( max(\psi_l \dot{-} \frac{k'_l }{l}, 
\frac{k_l}{l} \dot{-} \psi_l ,l_{n,\phi} \dot{-} (\phi + \sum_{i\le m} 2^{-(2n+i+1)} ), 
\psi \dot{-} \frac{k_{m+1} +2}{m+1})). 
$$ 
Since $T$ is decidable we compute 
$k'_{m+1} <m$ so that 
$\frac{k'_{m+1}}{m+1} \le (\theta' )^{\circ} \le \frac{k'_{m+1}+2}{m+1}$.  
Then the value of $\psi$ under $T_{n,\phi ,m}$ 
together with $\psi \le \frac{k_{m+1} +2}{m+1}$
is not greater than $\frac{k'_{m+1}+2}{m+1}$. 
This means that extending $T_{n,\phi ,m}$ by 
$\psi \le \frac{min(k_{m+1}, k'_{m+1})+2}{m+1}$ 
we preserve consistency of the theory.

If $0< k'_{m+1} <k_{m+1}$ 
we repeat this construction again. 
It is clear that finally we arrive at the situation 
when after such a repetition the number $k_{m+1}$ does not change.  
Then note that the value $\psi^{\circ}$ under 
the extension of $T$ by  
$\psi \le  \frac{k_{m+1}+2}{m+1} + 2^{-(2n+m+2)}$ 
and all statements of the form 
$$
\frac{k_l}{l} -2^{-(2n+m+2)}  \le \psi_l \le \frac{k'_l }{l} +2^{-(2n+m+2)}  
\mbox{ (for inequalities } \frac{k_l}{l} \le \psi_l \le \frac{k'_l}{l}
\mbox{  from } T_{n,\phi ,m} ) 
$$ 
satisfies  $\frac{k_{m+1}}{m+1}\le \psi^{\circ}$ 
(apply $\frac{k_{m+1}}{m+1} \le \theta^{\circ}$ with respect to $T$).  
We now define $T_{n, \phi ,m+1}$ 
as the set of so corrected statements of $T_{n,\phi ,m}$ 
together with the statement 
$$
\frac{k_{m+1}}{m+1}\le \psi \le \frac{k_{m+1}+2}{m+1}+ 2^{-(2n+m+2)} . 
$$ 
It is clear that the obtained extension is consistent with $T$. 

By the choice of a repeating enumeration
we see that for each sentence $\psi$ 
boundaries of $\psi$ at steps of our procedure 
form a Cauchy sequence. 
Thus $\psi$ has the same value in all models of 
$T_{n,\phi}$. 
Moreover the inequality $l_{n,\phi} \dot{-} \phi \le 0$ 
will be transformed into $l_{n,\phi} \dot{-} \phi \le 2^{-n}$. 
We see that Step 0 guarantees that $T$ coincides 
with the intersection  of all $T_{n,\phi}$.  

Note that after the $(m+1)$-th step 
we know that for every inequality $\psi'\le \delta$ 
from each $T_{n',\phi' ,m+1}\setminus T$ 
the upper boundary of $\psi'$ in the final 
$T_{n' ,\phi'}$ cannot exceed $\delta + \frac{1}{2^{m}}$. 
In particular all ineqalities of this kind can be included 
into an enumeration of axioms of $T_{n',\phi'}$ at this step.   
Thus we see that by the effectiveness of our procedure 
the family $\{ T_{n,\phi} \}$ is effective.  
$\Box$ 

\bigskip 

In order to have a method for 
proving undecidability of continuous 
theories we now discuss 
interpretability of first order structures 
in continuous ones.  

Let $L_0 =\langle P_1 ,...,P_m  \rangle$ be a finite 
relational signature.  
Let $\mathcal{K}_0$ be a class 
of finite first-order $L_0$-structures. 
Let $\mathcal{K}$ be a class of continuous 
$L$-structures, where $L$ is as above. 
We say that $\mathcal{K}_0$ is {\bf relatively 
interpretable} in $\mathcal{K}$ if there is a finite constant 
extension $L(\bar{a}) = L\cup \{ a_1 ,...,a_r \}$, 
a constant expansion $\mathcal{K} (\bar{a})$ of $\mathcal{K}$ 
and there  are continuous $L$-formulas 
$$ 
\phi^{-} (\bar{x}, \bar{y}) \mbox{ , } \phi^{+} (\bar{x}, \bar{y}) 
\mbox{ , } \theta^{-} (\bar{x}, \bar{y}_1 , \bar{y}_2) 
\mbox{ , } \theta^{+} (\bar{x}, \bar{y}_1 , \bar{y}_2)  \mbox{ and } 
$$ 
$$ 
\psi^{-}_1 (\bar{x}, \bar{y}_1 , \bar{y}_2 ,...,\bar{y}_{l_1} ) 
\mbox{ , } \psi^{+}_1 (\bar{x}, \bar{y}_1 , \bar{y}_2 ,...,,\bar{y}_{l_1} ) \mbox{  , ..., } 
\psi^{-}_m (\bar{x}, \bar{y}_1 , \bar{y}_2 ,...,\bar{y}_{l_m} ) 
\mbox{ , } \psi^{+}_m (\bar{x}, \bar{y}_1 , \bar{y}_2 ,...,,\bar{y}_{l_m} ) , 
$$ 
$$ 
\mbox{ with  } 
|\bar{y} |= |\bar{y}_1 | = |\bar{y}_2 | = ... |\bar{y}_{l_j}| = ... =|\bar{y}_{l_m} | 
\mbox{ , such  that: } 
$$ 
(i) the $L$-reduct of $\mathcal{K}(\bar{a})$ coincides with $\mathcal{K}$;\\ 
(ii) the conditions $\phi^{-} (\bar{a}, \bar{y})\le 0$ and 
$\phi^{+} (\bar{a}, \bar{y}) > 0$ are equivalent  in any $M\in \mathcal{K}(\bar{a})$
and  the condition $\theta^{-} (\bar{a}, \bar{y}_1, \bar{y}_2)\le 0$  
defines an equivalence relation on the zero-set of $\phi^{-} (\bar{a} ,\bar{y})$ 
(on tuples of the corresponding power $M^s$ with $s = |\bar{y}_1 |$),  
so that the values of any 
$\psi^{\varepsilon}_i (\bar{a},\bar{y}_1 , \bar{y}_2 ,...,\bar{y}_{l_i} )$ 
are invariant under this equivalence relation; \\ 
(iii)  the $(+)$-conditions below are equivalent to $(-)$-ones in $\mathcal{K}(\bar{a})$ :  
$$ 
\theta^{-} (\bar{a}, \bar{y}_1, \bar{y}_2)\le 0 \mbox{ , }  
\theta^{+} (\bar{a}, \bar{y}_1 ,\bar{y}_2) > 0 \mbox{ , } 
\psi^{-}_1 (\bar{a}, \bar{y}_1 , \bar{y}_2 ,...,\bar{y}_{l_1} ) \le 0 
\mbox{ , } 
$$ 
$$ 
\psi^{+}_1 (\bar{a}, \bar{y}_1 , \bar{y}_2 ,...,,\bar{y}_{l_1} )>0 \mbox{  , ..., } 
\psi^{-}_m (\bar{a}, \bar{y}_1 , \bar{y}_2 ,...,\bar{y}_{l_m} ) \le 0
\mbox{ , } 
\psi^{+}_m (\bar{a}, \bar{y}_1 , \bar{y}_2 ,...,,\bar{y}_{l_m} )>0 ;  
$$ 
(iv) for any $M \in \mathcal{K}(\bar{a})$ the conditions of (iii)  define 
an $L_0$-structure from $\mathcal{K}_0$ on the $\theta$-quotient   
of the zero-set  of $\phi^{-} (\bar{a} ,\bar{y})$  and 
any structure of $\mathcal{K}_0$ can be so realised.  

\begin{thm} \label{unde}  
Under circumstances above assume that $Th(\mathcal{K}_0 )$ 
is undecibable. 
Then  the continuous theory  $Th(\mathcal{K}(\bar{a}))$ 
is not a computable set. 
\end{thm} 

{\em Proof.} 
The proof is straightforward. 
To each formula $\psi$ of the theory of $\mathcal{K}_0$  so that 
the quantifier-free part is in the disjunctive normal form we associate 
the appropriately rewritten continuous formula $\psi^- (\bar{a})$ 
(with appropriate free variables).  
In particular atomic formulas are written by $(-)$-conditions 
above, but negations of atomic formulas 
appear in the form of 
$$
\psi^{+}_i (\bar{a}, \bar{y}_1 , \bar{y}_2 ,...,,\bar{y}_{l_i} )\le 0 . 
$$  
Condition (ii) and the condition that  the $\theta$-quotient   
of the zero-set  of $\phi^{-} (\bar{a} ,\bar{y})$ 
is always finite, allow us to use standard quantifiers 
in such $\psi^- (\bar{a})$: the quantifier $\forall$ 
is written as $sup$ but $\exists$ is written as $inf$.  
Note that if $\psi'$ is equivalent to $\neg \psi$  
then  $(\psi' )^{-} (\bar{a})\le 0$ is equivalent 
to  $\psi ^{-} (\bar{a})> 0$ for tuples from 
the zero-set of $\phi^{-} (\bar{a}, \bar{y})$. 

This construction reduces the decision problem 
for $Th(\mathcal{K}_0 )$ to  computability of 
$Th(\mathcal{K}(\bar{a}))$.  
$\Box$

\bigskip 

It is worth noting that in the classical first-order logic 
the situation of this theorem usually has much stronger consequences. 
For example Theorem 5.1.2 of \cite{ershov} in a slightly modified 
setting (and removing the assumption that $\mathcal{K}_0$ consists of 
finite structures) states that hereditary undeciadability of 
$Th(\mathcal{K}_0)$ can be lifted to $Th(\mathcal{K})$.   
The 'positiveness' of the continuous logic does not allow so strong statements.  
  
As we already know the statement of Theorem \ref{unde} 
does not imply that $Th(\mathcal{K}(\bar{a}))$ is undeciadable. 
To prove the undeciadability theorems from Section 4 
we will use an additional tool.  
It will be applied in a combination with the method of Theorem \ref{unde}. 


\section{Decidability/undecidability of theories of qubit spaces} 

We start this section with an undecidability 
result of some classes of constant expansions of qubit spaces. 
Then we apply the idea of the proof 
to a more interesting example of a class of dynamical 
qubit spaces.

\begin{thm}  \label{constants} 
There is a class of qubit spaces expanded by four constants, 
i.e. structures of the form
$$  
(Q, qu , \{ B_n\}_{n\in \omega} ,0,\{ I_{mn} \}_{m<n} ,
\{ \lambda_c \}_{c\in\mathbb{Q}[i]}, +,-,\langle \rangle_{Re} , \langle \rangle_{Im}, a_1 ,a_2 , b_1 ,b_2  )  
$$
which is distinguished in the class of all qubit spaces 
with $\langle b_1 ,b_2 \rangle \not= 0$ by a continuous 
statement and which has undecidable continuous theory. 
\end{thm} 

{\em Proof.} 
Consider the following formula: 
$$ 
\psi (x ,y_1 ,y_2 ) = |\langle qu(y_1 ),x\rangle - \langle qu(y_2 ),x \rangle | ,  
$$ 
where $y_1 ,y_2$ are variables of the sort $Q$ and $x$ is of the sort $B_1$. 

In any qubit space any equivalence relation on $Q$ can be realised by 
$\psi (a , y_1 ,y_2 )\le 0$ for appropriately chosen $a\in B_1$:   
define $a$ to be a linear combination of $qu (q_l )$, so that for 
equivalent $q_j$ and $q_k$ the coefficients of $qu (q_j )$ and $qu (q_k )$ 
in $a$ are the same. 

Let us introduce the following formula: 
$$ 
\psi^c (x , z_1 ,z_2 ,y_1 ,y_2 ) = 
| \langle z_1 ,z_2 \rangle |\dot{-} |\langle qu(y_1 ),x\rangle - \langle qu(y_2 ),x \rangle | , 
$$ 
where $y_1 ,y_2$ are variables of the sort $Q$ and $x,z_1 ,z_2$ are of $B_1$.

If $a$ defines an equivalence relation on $Q$ as above then there are 
$b_1 ,b_2 \in B_1$ with $\langle b_1 ,b_2 \rangle \not=0$ and 
$$ 
sup_{y_1 ,y_2} min(\psi (a, y_1 ,y_2 ), \psi^c (a,b_1 ,b_2 , y_1 , y_2 ))\le 0
$$ 
$$
sup_{y_1 ,y_2} (| \langle b_1 ,b_2 \rangle |\dot{-} (\psi (a, y_1 ,y_2 )+ \psi^c (a,b_1 ,b_2 , y_1 , y_2 )))\le 0.
$$
We see that the formula 
$\psi^c (a,b_1 ,b_2, y_1 , y_2 )$ can be 
interpreted as the complement of the equivalence relation 
defined by $\psi (a, y_1 , y_2 )$ in the class of these qubit spaces. 

This allows us to define interpretability of the first-order theory of 
finite structures of two equivalence relations  
(which is undecidable by Proposition 5.1.7 from \cite{ershov}) 
in the class, say $\mathcal{K}$, of qubit spaces extended by constants 
$a_1 ,a_2 ,b_1 ,b_2$ where $b_1 ,b_2$ satisfy the statements 
above for both $a_1$ and $a_2$ instead of $a$. 

In fact we axiomatise $\mathcal{K}$ in the class of all qubit spaces 
with $\langle b_1 , b_2 \rangle \not= 0$ by the following continuous statements: 
$$ 
sup_{y_1 ,y_2} min(\psi (a_i , y_1 ,y_2 ), 
| \langle b_1 ,b_2 \rangle |\dot{-} |\langle qu(y_1 ),a_i \rangle - \langle qu(y_2 ),a_i \rangle |) \le 0 
\mbox{ , where } i=1,2.  
$$ 

In terms of  Theorem \ref{unde} the formulas $\phi^+ , \phi^- , \theta^+ , \theta^-$ 
become degenerate: $\phi^-$ can be taken as $d(y,y)$ for 
the (descrete) sort $Q$, then 
$$
\phi^+ (y) = 1 \dot{-} d(y,y) \mbox{ , } 
\theta^- (y_1 ,y_2 ) = d(y_1 , y_2 ) \mbox{ , } 
\theta^+ (y_1 ,y_2 ) = 1 - d(y_1 ,y_2 ). 
$$ 
Formulas   $\psi (a_1 , y_1 , y_2 )$ and $\psi^c (a_1 ,b_1 ,b_2, y_1 , y_2 )$ 
play the role of $\psi^{-}_1$  and $\psi^+_1$. 
Then $\psi (a_2 , y_1 , y_2 )$ and $\psi^c (a_2 ,b_1 ,b_2, y_1 , y_2 )$
play the role of  $\psi^{-}_2$  and $\psi^+_2$. 

To each formula $\rho (\bar{y})$ of the theory of two equivalence relations 
so that the quantifier-free part is in the disjunctive normal form we associate 
the appropriately rewritten continuous formula 
$\rho^* (a_1 ,a_2 ,b_1 ,b_2 , \bar{y})$ 
(where we use $min$ and $max$ instead of $\vee$ and $\wedge$). 
Since the free variables of the latter $\bar{y}$ are of the sort $Q$,  
when $\rho$ is quantifier-free,  the values 
$\rho^* (a_1 ,a_2 ,b_1 ,b_2 , \bar{c})$ belong to 
$\{ 0\} \cup [ | \langle b_1 ,b_2 \rangle |, 1]$.  
It is easy to see that in structures of $\mathcal{K}$ 
the same property holds for any 
formula $\rho (\bar{y})$. 

This obviously implies that when 
$\rho$ is a sentence,  the sentence 
$$ 
min (| \langle b_1 ,b_2 \rangle |,\rho^*(a_1 ,a_2 ,b_1 ,b_2 ))
$$   
has the following property:  
\begin{quote}  
$\rho$ is satisfied in all finite models of two equivalence relations 
if and only if all structures of $\mathcal{K}$  satisfy 
$min (| \langle b_1 ,b_2 \rangle |,\rho^*(\bar{a},\bar{b})) =0$.  
\end{quote} 
If the theory of $\mathcal{K}$ was decidable, it would 
define the real number $| \langle b_1 ,b_2 \rangle |^{o}$ which 
would be computable and $\not= 0$. 
In particular we could compute a rational $r>0$, so that 
$min (| \langle b_1 ,b_2 \rangle |,\rho^*(\bar{a},\bar{b})) =0$ 
is equivalent to $min (| \langle b_1 ,b_2 \rangle |,\rho^*(\bar{a},\bar{b}))\le r$.  
By Theorem \ref{unde} this is impossible. 
$\Box$ 

\bigskip

\begin{thm} 
There is a class of dynamical qubit spaces  
in the signature 
$$  
(Q, qu , \{ B_n\}_{n\in \omega} ,0,\{ I_{mn} \}_{m<n} ,
\{ \lambda_c \}_{c\in\mathbb{Q}[i]}, +,-,\langle \rangle_{Re} , 
\langle \rangle_{Im}, U_1 ,U_2 , U_3 ,U_4 ,U_5  )  
$$
which is distinguished in the class of all qubit spaces with 
$$
sup_v d( U_3 (v) ,v )\not= 0
$$ 
by a continuous statement and 
which has undecidable continuous theory. 
\end{thm} 

{\em Proof.} 
We will use the construction of Theorem \ref{constants} 
with some necessary changes. 
For example we replace the value  
$| \langle b_1 ,b_2 \rangle |$ from  
that theorem  by  $sup_{v} d(U_3 (v ),v )$. 
The  constants $a_1$ and $a_2$  will appear as 
the normed vectors fixed by $U_1$ and $U_2$ 
respectively. 
Although we choose $U_i$, $i=1,2$, so that 
the subspace of fixed vectors of $U_i$ 
coincides with $\mathbb{C} a_i$, 
we cannot define these constants by a continuous formula. 
This is why some additional values will be 
used in the proof.  
The values  $sup_{v} d(U_3 (v ),v )$  
and $sup_{v} d(U_4 (v ),v )$ will appear 
in 'fuzzy' versions of formulas from 
the proof of Theorem \ref{constants}. 

Let: 
$$ 
\psi_i (y_1 ,y_2 ) = sup_{u} min ( 
sup_{v_1}(d(U_{3} (v_1 ),v_1 )) \dot{-} max (d(U_i (u ),u ),|1-\parallel u \parallel |),
$$ 
$$    
(|\langle qu(y_1 ),u\rangle - \langle qu(y_2 ),u \rangle | \dot{-} 
sup_{v_2}d(U_4 (v_2 ),v_2 ))),  
$$ 
where $y_1 ,y_2$ are variables of the sort $Q$, $i\in \{ 1,2\}$  
and $u , v_1 ,v_2$ are of the sort $B_1$. 

To see that in any qubit space any equivalence 
relation on $Q$ can be realised by $\psi_i (y_1 ,y_2 )\le 0$ 
let us define $a_1$ (the case of $a_2$ is similar) 
to be a linear combination of 
$qu (q_l )$ of length 1, so that for equivalent 
$q_j$ and $q_k$ the coefficients of $qu (q_j )$ 
and $qu (q_k )$ in $a_1$ are the same. 
We also fix a rational number $r$ (for both $a_1$ and $a_2$)  
so that for non-equivalent $q_j$ and $q_k$ the coefficients 
of $qu (q_j )$ and $qu (q_k )$ in $a_1$ 
are distant by $> r$.   
Note that any $e^{i\phi} a_1$ has the same properties 
with respect to elements of $qu(Q)$. 

We extend $a_1$ to an orthonormal basis of 
$\mathbb{H}$ and define $U_1$ to be an unitary 
operator having these vectors as eigenvectors 
so that $\mathbb{C} a_1$ is the subspace of fixed points. 
The remaining eigenvalues are chosen so that 
the corresponding eigenvectors are taken by $U_i$ 
at the distance $\ge 1/10$. 

Now we can take the operators $U_3$ and $U_4$  
so close to $Id$ (with respect to the operator norm) 
that the formula $\psi_1$ indeed realises 
the equivance relation we consider.  
Choosing $U_4$ we demand that 
$sup_{v_2}(d(U_4 (v_2 ),v_2 ))$ 
is much less than $r$. 
Having $U_4$ we can find $U_3$. 
We will assume that $sup_{v}d(U_{3} (v),v )>0$ 
in our structures. 

Let us now introduce $U_5$ with 
$r =  sup_{v}(d(U_{5} (v ),v)$ 
and consider the following formulas for $i = 1,2$: 
$$ 
\psi^c_i (y_1 ,y_2 ) =  sup_{u} min ( 
sup_{v_1}(d(U_{3} (v_1 ),v_1 )) \dot{-} max (d(U_i (u ),u ),|1-\parallel u \parallel |),
$$ 
$$ 
(sup_{v_2}d (U_5 (v_2 ),v_2 )\dot{-} |\langle qu(y_1 ),u\rangle - \langle qu(y_2 ),u \rangle | ) ), 
$$ 
where $y_1 ,y_2$ are variables of the sort $Q$ and $u,v_1 ,v_2$ are of $B_1$.
If necessary we may correct $U_3$ making $sup_{v}d(U_{3} (v),v)$ smaller  
so that the following statements hold. 
$$ 
sup_{y_1 ,y_2} min(\psi_i (y_1 ,y_2 ), \psi^c_i ( y_1 , y_2 ))\le 0 ,
$$ 
$$
sup_{y_1 ,y_2} ( sup_{v}d(U_{3} (v),v )\dot{-} (\psi_i ( y_1 ,y_2 )+ \psi^c_i ( y_1 , y_2 )))\le 0.
$$
As before the formula $\psi^c_i (y_1 , y_2 )$ will be 
interpreted as the complement of the equivalence relation 
defined by $\psi_i (y_1 , y_2 )$ in the class of these qubit spaces. 

This allows us to define interpretability of the (undecidable) 
first-order theory of finite structures of two equivalence relations in 
the class, say $\mathcal{K}$, of dynamical qubit spaces with respect to 
operators  $U_1 ,U_2 ,U_3 ,U_4 , U_5$. 
The formulas 
$\phi^+ , \phi^- , \theta^+ , \theta^-$ 
(see  Theorem \ref{unde} ) are taken as in 
Theorem \ref{constants} (i.e. $\phi^- (y) = d(y,y)$  
and $\theta^- (y_1 ,y_2 ) = d(y_1 , y_2 )$).  
Formulas   $\psi_i (y_1 , y_2 )$ and $\psi^c _i ( y_1 , y_2 )$ 
play the role of $\psi^{-}_i$  and $\psi^+_i$ for  $i=1, 2$. 

To each formula $\rho (\bar{y})$ of the theory of two equivalence relations 
so that the quantifier-free part is in the disjunctive normal form we associate 
the appropriately rewritten continuous formula 
$\rho^* ( \bar{y})$. 
Since the free variables of the latter $\bar{y}$ are of the sort $Q$,  
when $\rho$ is quantifier-free,  the values 
$\rho^* ( \bar{c})$ belong to 
$\{ 0\} \cup [ sup_{v}d(U_{3} (v),v ), 1]$.  
Thus we see that in structures of $\mathcal{K}$ 
the same property holds for any 
formula $\rho (\bar{y})$. 

This obviously implies that when 
$\rho$ is a sentence,  the sentence 
$$ 
min (sup_{v}d(U_{3} (v),v ),\rho^*)
$$   
has the following property:  
\begin{quote}  
$\rho$ is satisfied in all finite models of two equivalence relations 
if and only if all structures of $\mathcal{K}$  satisfy 
$min (sup_{v}d(U_{3} (v),v ),\rho^*) =0$.  
\end{quote} 
If the theory of $\mathcal{K}$ was decidable, it would 
define the real number $sup_{v}d(U_{3} (v),v )^{o}$ which 
would be computable and $\not= 0$. 
In particular we could compute a rational $r'>0$, so that 
$min (sup_{v}d(U_{3} (v),v ),\rho^* ) =0$ 
is equivalent to $min (sup_{v}d(U_{3} (v),v ),\rho^* )\le r'$.  
By Theorem \ref{unde} this is impossible. 
$\Box$ 

\bigskip


Let us now restrict the dimension 
of qubit spaces, say by $2^n$. 
It is natural to expect that 
then the theory of (dynamical) 
qubit spaces becomes decidable. 
In classical model theory this corresponds 
to the situation of a theory of structures 
of a fixed finite size. 
 
On the other hand  since the structures 
are of infinite language it is not very difficult 
to find such a structure 
with undecidable continuous theory. 
For example one can take a dynamical qubit space 
with one additional operator which takes a basic vector, 
say $qu (q_0 )$, to some linear combination 
$r \cdot qu(q_0 ) + \sqrt{1 - r^2} \cdot qu (q_1 )$ 
where $r$ is a non-computable real number. 

Let us fix a signature 
$$
(Q, qu , \{ B_n\}_{n\in \omega} ,0,\{ I_{mn} \}_{m<n} ,
\{ \lambda_c \}_{c\in\mathbb{Q}[i]}, +,-,\langle \rangle_{Re} , \langle \rangle_{Im}, U_1 ,...,U_t  ),
$$
where as before we assume that $U_j$, $1\le j \le t$,  
are  symbols of unitary operators of $\mathbb{H}$ 
which are defined only on $B_1$. 
Using Theorem \ref{Ersov} we will 
prove that the theory of $2^n$-qubit 
spaces in this language is decidable. 

Let us enumerate all $2^n$-dimensional 
unitary matrices of computable complex numbers 
(i.e. their real and imaginary parts are defined 
by computable sequences). 
This induces an enumeration $Axm_j$, $j\in \omega$, 
of systems of axioms of complete continuous theories  
$T_j$ of dynamical $n$-qubit spaces. 
Each $Axm_j$  consists of the standard 
axioms of $n$-qubit spaces, the axioms stating that each 
$U_s$ is a unitary operator and the axioms 
describing the matrix of $U_s$ in the 
basis $qu (Q)$: 
$$ 
inf_{q_1 ...q_N} \parallel U_s (qu (q_l ))-\Sigma \lambda_{c_j}(qu (q_j )) \parallel \le \varepsilon_l 
\mbox{ , where } 
\varepsilon\in \mathbb{Q}, c_j \in \mathbb{Q}[i] . 
$$ 
Using Lemma \ref{comp} it is easy to see that the enumeration $Axm_j$, $j\in \omega$, 
gives an effective indexation of complete continuous theories  
$T_i$ of dynamical $n$-qubit spaces in the sense of Section 3.    
The statement that the relation 
$\{ (\theta ,j): \theta$ is a statement so that $T_j \vdash \theta \}$ 
is computably enumerable follows from   
the fact that ths relation coincides with 
$\{ (\theta ,j): \theta$ is a statement so that $Axm_j \vdash \theta \}$.   

\begin{thm} 
The theory of all dynamical $n$-qubit spaces 
coincides with the intersection $\bigcap T_j$.  \\ 
The theory of all dynamical $n$-qubit spaces is decidable. 
\end{thm} 

{\em Proof.} 
By Theorem \ref{Ersov} the second statement of 
the theorem follows from the first one. 
Thus we only have to show that 
for any rational $\delta$, any dynamical $n$-qubit space  
$$
(Q, qu ,\{ B_n\}_{n\in \omega} ,0,\{ I_{mn} \}_{m<n} ,
\{ \lambda_c \}_{c\in\mathbb{Q}[i]}, +,-,\langle \rangle_{Re} , \langle \rangle_{Im}, U_1 ,...,U_t  ),
$$ 
and any continuous sentence $\theta ( U_1 ,...,U_t  )$ over this 
structure there are operators $\tilde{U}_1 ,...,\tilde{U}_t$ 
defined by matrices over $\mathbb{Q}[i]$, so that  
$$ 
| \theta ( U_1 ,...,U_t  ) - \theta ( \tilde{U}_1 ,...,\tilde{U}_t  )| \le \delta . 
$$ 
Indeed this shows that when some $\theta ( U_1 ,...,U_t  ) \le \varepsilon$ 
does not belong to $T$, then it does not belong 
to some $T_j$. 

Since any continuous formula defines a uniformly 
continuous function and the ball $B_1$ is compact 
it suffices to take $\tilde{U}_1 ,...,\tilde{U}_t$ so that 
they sufficiently approximate $U_1 ,...,U_t$.  
This is a folklore fact. 
On the other hand it is a curious place where 
the following fact from quantum computations can be applied. 

Let $CNOT$ be a 2-qubit linear operator defined by  
$$ 
CNOT: |00\rangle \rightarrow |00\rangle \mbox{ , } |01\rangle \rightarrow |01\rangle 
\mbox{ , } 
|10\rangle \rightarrow |11\rangle \mbox{ , } |11\rangle \rightarrow |10\rangle   
. 
$$ 
The {\em Toffoli gate} is a 3-qubit linear operator defined on basic vectors by 
$$
\Lambda(CNOT): |\varepsilon_1 \varepsilon_2 \varepsilon_3 \rangle \rightarrow 
|\varepsilon_1 \varepsilon_2 (\varepsilon_3 \oplus \varepsilon_1 \cdot \varepsilon_2 )\rangle .  
$$ 
Let $\mathcal{B}$ be a 2-dimensional space over $\mathbb{C}$. 
It is well-known that   
\begin{quote} 
(a) Any unitary transformation of $(\mathcal{B} )^{\otimes k}$ is a product of  
1-qubit unitary transformations and  2-qubit copies of CNOT  at appropriate registers. 

(b)  The operators of the basis  
$$ 
\mathcal{Q}  = \{  K = {{ 1 \mbox{ } 0}\choose{0 \mbox{ } i}} , CNOT , \Lambda (CNOT), 
 \mbox{Hadamar's } H= \frac{1}{\sqrt{2}}{{ 1 \mbox{ } 1}\choose{1 \mbox{ } -1}} \} 
$$ 
generate a dense subgroup of  
$\mathbb{U}(\mathcal{B}^{\otimes 3})/\mathbb{U}(1)$ under the operator norm. $\Box$ 
\end{quote}

\bigskip

Institute of Mathematics, University of Wroc{\l}aw, \parskip0pt

pl.Grunwaldzki 2/4, 50-384 Wroc{\l}aw, Poland \parskip0pt

E-mail: ivanov@math.uni.wroc.pl


\begin{thebibliography}{99}
\bibitem{BYBHU} I.Ben Yaacov, A.Berenstein, W.Henson and A.Usvyatsov, 
{\em Model theory for metric structures}. In: Model theory with 
Applications to Algebra and Analysis, v.2 (Z.Chatzidakis, H.D.Macpherson, 
A.Pillay and A.Wilkie, eds.), London Math. Soc. Lecture Notes, v.350, 
pp. 315 - 427, Cambridge University Press, 2008.   
\bibitem {BYP} I.Ben Yaacov and A.R.Pedersen, 
{\em A proof of completeness for 
continuous first-order logic}, J.Symbolic Logic, 75(2010), 168 - 190. 
\bibitem{BYU} I.Ben Yaacov and A.Usvyatsov, {\em Continuous first order logic 
and local stability}, Trans Amer. Math. Soc., 362(2010), 5213 - 5259.  
\bibitem{CCGP} G.Cattaneo, M.L.Dalla Chiara, R.Giuntini and F.Paoli, 
{\em Quantum logic and nonclassical logics}. 
In: Handbook of Quantum Logic and Quantum Structures, Quantum Logic.  
(K.Engesser, D.M.Gabbay, D.Lehmann, eds.), 
pp. 127 - 226, Elsevier B.V., 2009. 
\bibitem{DJK} H.Derksen, E.Jeandel and P.Koiran, 
{\em Quantum automata and algebraic groups}, 
J. of Symbolic Computations 39(2005), 357 - 371. 
\bibitem{DGP} F.Didehvar, K.Ghasemloo and M.Pourmahdian, 
{\em Effectiveness in RPL, with applications to continuous logic}, 
Ann Pure Appl. Logic, 161(2010), 789 - 799. 
\bibitem{EH} A.Edalat and R.Heckmann, 
{\em A computational model for metric spaces}, 
Theoretical Computer Science, 193(1998), 53 - 73.
\bibitem{ES} A.Edalat and Ph.S\"{u}nderhauf, 
{\em A domain theoretic approach to computability on the real line}, 
Theoretical Computer Science, 210(1999), 73 - 98.  
\bibitem{ershov} Yu.L.Ershov, {\em Decision Problems and Constructive Models}. Nauka, Moscow, 1980. 
\bibitem{FHS} I.Farah, B.Hart and D.Sherman, {\em Model theory of operator algebras I:stability}, 
Bull. London. Math. Soc. 45(2013), no.4, 825 -838. 
\bibitem{GH} I.Goldbring and B.Hart, {\em A computability-theoretic reformulation of 
the Connes Embedding Problem}, arXiv: 1308.2638.  
\bibitem{gruska} J.Gruska, {\em Quantum Computing}. McGraw-Hill International (UK) Limited, London, 1999. 
\bibitem{MC} C.Moore and J.P.Crutchfield, {\em Quantum automata and quantum grammars}, 
Theoretical Computer Sci. 237(2000), 275 - 306.  
\bibitem{pykacz} J.Pykacz, {\em Unification of two approaches to quantum logic: 
every Bikhoff-von Neumann quantum logic is a partial infinite-valued {\L}ukasiewicz logic}, 
Studia Logica, 95 (2010), 5 - 20.  
\end{thebibliography}
\end{document}